\documentclass [oneside,a4paper,10pt] {article}
\usepackage{latexsym,epsfig,graphpap,graphics,amsmath,amssymb,}
\usepackage{hyperref}
\begin{document}

\begin{center}
{\bf{A NOTE ON QUASI UMBILICAL HYPERSURFACE OF A SASAKIAN MANIFOLD WITH $(\phi, g, u, v, \lambda)-$ STRUCTURE}}\\

\vspace{.3cm} \noindent{\bf{Sachin Kumar Srivastava and Alok Kumar Srivastava}}\\

\end{center}

\begin{abstract}
In this paper we have studied the properties of Quasi umbilical hypersurface $M$ of a Sasakian manifold $\tilde M$ with $(\phi, g, u, v, \lambda)-$structure and established the relation for $M$ to be cylindrical.
\end{abstract} 

\vspace{.3cm}\noindent{\bf{2000 Mathematics Subject Classification :}}~ 53D10, 53C25, 53C21.

\vspace{.5cm}
\noindent{\bf{1.~~Introduction~:}}~Let $\tilde M$ be a $(2n + 1)-$dimensional Sasakian manifold with a tensor field $\phi$ of type (1, 1),  a fundamental vector field $\xi$ and $1-form \, \eta$ such that\\
\nolinebreak
(1.1)\hspace{1in}$\eta(\xi) = 1$\\
(1.2)\hspace{1in}$\tilde \phi^2 = - I + \eta \otimes \xi $\\
\vspace{.2cm}
where $I$ denotes the identity transformation.\\
(1.3)(a)\hspace{.8in}
$\eta \, o\,\phi = 0\qquad (b) \qquad \phi \xi = 0 \qquad (c) \qquad rank(\tilde \phi) = 2n$\\
\vspace{.2cm}
If $\tilde M$ admits a Riemannian metric $\tilde g$, such that\\
\vspace{.2cm}
(1.4)\hspace{1in}$\tilde g(\tilde\phi X,\tilde\phi Y) = \tilde g(X, Y) - \eta(X) \eta(Y)$\\
\vspace{.2cm}
(1.5)\hspace{1in}$ \tilde g(X,\xi)=\eta(X)$\\
\vspace{.2cm}
then $\tilde M$ is said to admit a $(\tilde\phi,\xi,\eta,\tilde g)-$ structure called contact metric strucure.\\
If moreover,\\
(1.6)\hspace{1in}$(\tilde\nabla_X \tilde\phi )\, Y =\tilde g (X, \, Y) \xi - \eta (Y) X $\\
 \hspace{.3in} and\\
(1.7)\hspace{1in}$\tilde\nabla_X\xi=-\tilde\phi X$\\
where $\tilde\nabla$ denotes the Riemannian connection of the Riemannian metric g, then $(\tilde M,\tilde\phi,\xi,\eta,\tilde g)$ is called a Sasakian manifold ( see \cite{DEB285}).
\vspace{3mm}\parindent=8mm
If we define ${}^\prime F(X, Y) = g (\phi X, Y)$, then in addition to above relation we find\\
(1.8)\hspace{1in}${}^\prime F(X, Y) + {}^\prime F(Y, X) = 0 $\\
(1.9)\hspace{1in}${}^\prime F(X, \phi Y) = {}^\prime F(Y, \phi X)$\\
(1.10)\hspace{.92in}${}^\prime F(\phi X, \phi Y) = {}^\prime F(X, Y)$\\
\vspace{.5cm}
\noindent{\bf{2.~~Hypersurface of a Sasakian manifold with $(\phi, g, u, v, \lambda)-$structure~:}}~\\
\vspace{.1cm}
\hspace{.5cm}Let us consider a $2n-$dimensional manifold $M$ embedded in $\tilde M$ with embedding $b : M \rightarrow \tilde M$. The map $b$ induces a linear transformation map $B$ (called Jacobian map), $B : T_p \rightarrow T_{b_p}$.
\vspace{.1cm}
Let an affine normal $N$ of $M$ is in such a way that $\tilde\phi N$ is always tangent to the hypersurface and satisfying the linear transformations

\vspace{.25cm}\noindent
(2.1) \hspace{1in} $\tilde\phi BX = B \phi X + u(X) N$\\
(2.2) \hspace{1in} $\tilde\phi N = - BU$\\
(2.3) \hspace{1in} $\xi = BV + \lambda N$\\
(2.4) \hspace{1in} $\eta(BX) = v(X)$

\vspace{.5cm}\noindent
where $\phi$ is a (1, 1) type tensor; $U,\, V$ are vector fields; $u, v$ are $1- form$ and $\lambda$ is a $C^\infty -$ function. If $u \ne 0$, $M$, is called a noninvariant hypersurface of $\tilde M$ (see\cite{SIG25} and \cite{DNSK1611}).

\vspace{.3cm}\parindent=8mm
Operating (2.1), (2.2), (2.3) and (2.4)  by $\tilde\phi$ and using (1.1), (1.2) and (1.3) and taking tangent normal parts separately, we get the following induced structure on $M$:

\vspace{.3cm}\noindent
(2.5)(a) \hspace{.96in} $\phi^2 X = -  X + u (X) U + v (X) V$

\vspace{.3cm}\parindent=8mm
(b)\hspace{1in} $u(\phi X) = \lambda v (X), \qquad v(\phi X) =  -\, \eta (N) \, u(X)$

\vspace{.3cm}\parindent=8mm
(c)\hspace{1in} $\phi U = - \, \eta(N) \, V , \qquad \phi V = \lambda U$

\vspace{.3cm}\parindent=8mm
(d)\hspace{1in} $u(U) = 1 -  \lambda \eta(N), \qquad u(V) = 0$

\vspace{.3cm}\parindent=8mm
(e)\hspace{1in} $v(U) = 0, \qquad v(V) = 1 -   \lambda \eta(N)$

\vspace{.3cm}\noindent
and from (1.4) and (1.5), we get the induced metric $g$ on $M$, i.e.,

\vspace{.3cm}
\noindent
(2.6) \hspace{1in} $g(\phi X, \phi Y) = g (X, Y) - u(X) u(Y) - v(X) v(Y)$\\
(2.7) \hspace{1in} $g (U, X) = u (X), \qquad g (V, X) = v (X).$

\vspace{.3cm}
\parindent=8mm
If we consider $\eta(N) = \lambda$, we get the following structures on $M$.

\vspace{.3cm}
\noindent
(2.8)(a)\hspace{1.1in} $\phi^2 = -  I +  u \otimes U + v \otimes V$

\vspace{.3cm}
\parindent=8mm
(b) \hspace{1in} $\phi U = - \, \lambda V, \qquad \phi V = \lambda U$

\vspace{.3cm}
\parindent=8mm
(c) \hspace{1in} $u \, \circ \phi = \lambda v, \qquad v \, \circ \, \phi =  - \lambda u$

\vspace{.3cm}
\parindent=8mm
(d) \hspace{1in} $u( U) = 1 - \lambda^2 , \qquad u(V) = 0$

\vspace{.3cm}
\parindent=8mm
(e) \hspace{1in} $v (U) = 0,  \qquad v (V) = 1 - \lambda^2.$

\vspace{.3cm}
\parindent=8mm
A manifold $M$ with a metric $g$ satisfying (2.6), (2.7) and (2.8) is called manifold with $(\phi, g, u, v, \lambda )-$structure ( see\cite{DN799} and \cite{BBDN267}).
\vspace{.3cm}
\parindent=8mm
Let $\nabla$ be the induced connection on the hypersurface $M$ of the affine connection $\tilde\nabla$ of $\tilde M$.

\vspace{.3cm}
\parindent=8mm
Now using Gauss and Weingarten's equations

\vspace{.3cm}
\noindent
(2.9)\hspace{1in} $\tilde\nabla _{BX} BY = B\nabla_X Y + h (X, Y) N$

\vspace{.3cm}
\noindent
(2.10) \hspace{.9in} $\tilde\nabla_{BX} N = BHX + w (X) N,$ ~where~ $g (HY, Z) = h(Y, Z).$

\vspace{.3cm}
\parindent=8mm
Here $h$ and $H$ are the second fundamental tensors of type (0, 2) and (1, 1) and $w$ is a $1-form$. Now differentiating (2.1), (2.2), (2.3) and (2.4) covariantly and using (2.9), (2.10), (1.6) and reusing (2.1), (2.2), (2.3) and (2.4), we get

\vspace{.3cm}
\noindent
(2.11) \hspace{.9in} $(\nabla_Y \phi) (X) = v (X) Y - g (X, Y) V - h (X, Y) U - u (X) HY$

\vspace{.3cm}
\noindent
(2.12) \hspace{.9in} $(\nabla_Y u) (X) = - h(\phi X, Y) - u (X) w (Y) - \lambda g (X, Y)$

\vspace{.3cm}
\noindent
(2.13) \hspace{.9in} $(\nabla_Y v) (X) = g(\phi Y, X) + \lambda h (X, Y)$

\vspace{.3cm}
\noindent
(2.14) \hspace{.9in} $\nabla_Y U = w (Y) U - \phi HY - \lambda Y$

\vspace{.3cm}
\noindent
(2.15) \hspace{.9in} $\nabla_Y V = \phi Y + \lambda H Y$

\vspace{.3cm}
\noindent
(2.16) \hspace{.9in} $h(Y, V) = u(Y) - Y \lambda - \lambda w(Y)$

\vspace{.3cm}
\noindent
(2.17) \hspace{.9in} $h(Y, U) = - u(HY)$

\vspace{.3cm}
\parindent=8mm
Since $h (X, Y) = g (HX, Y)$, then from (1.5), and (2.17), we get

\vspace{.3cm}
\noindent
(2.18) \hspace{.9in} $h (Y, U) = 0 \Rightarrow HU = 0.$\\\\
\vspace{.1cm}
\noindent{\bf{3.~~Quasi-Umbilical hypersurface~:}}~If

\vspace{.3cm}
\noindent
(3.1) \hspace{1in} $h(X, Y) = \alpha g (X, Y) + \beta q (X) q(Y)$

\vspace{.3cm}
\noindent
where $\alpha, \beta$ are scalar functions, $q$ is $1-form$, then $M$ is called Quasi-umbilical hypersurface (see \cite{DEB285}) if  there exist a vector field $Q$ such that $g(Q, X) = q(X)$, where $g$ is the induced metric on $M$. If $\alpha = 0$, $\beta \ne 0$, then Quasi umbilical hypersurface $M$ is called {\it cylindrical} hypersurface (see \cite{DEB285}). If $\alpha \ne 0$, $\beta = 0$, then Qquasi-umbilical hypersurface $M$ is called {\it totally umbilical} and if $\alpha = 0, \beta = 0$, then Quasi umbilical hypersurface is {\it totally geodesic}.

\vspace{.3cm}
\parindent=8mm
Using (3.1) in (2.11), (2.12), (2.13), (2.14), (2.15), (2.16) and (2.17) we get\\\\
\vspace{.2cm}
(3.2)\hspace{.5in} $(\nabla_Y \phi) (X) = v (X) Y - g (X, Y) V - \{\alpha g (X, Y) + \beta q (X) q(Y)\} U $
$$\hspace*{2.2in} - u (X) \{\alpha Y + \beta q (Y) Q\}$$
\vspace{.2cm}
(3.3) \hspace{.5in} $(\nabla_Y u) (X) = - \{\alpha g (\phi X, Y) + \beta q (\phi X) q(Y)\}  - u (X) w (Y) - \lambda g (X, Y)$\\
\vspace{.2cm}
(3.4) \hspace{.5in} $(\nabla_Y v) (X) = g(\phi Y, X) + \lambda \{\alpha g (X, Y) + \beta q (X) q(Y)\}$\\
\vspace{.2cm}
(3.5) \hspace{.5in} $\nabla_Y U = w (Y) U - \{\alpha \phi Y + \beta q(Y) Q\} - \lambda Y$\\
\vspace{.2cm}
(3.6) \hspace{.5in} $\nabla_Y V = \phi Y + \lambda \{\alpha Y + \beta q(Y) Q\}$\\
\vspace{.2cm}
(3.7) \hspace{.5in} $h(Y, V) = \alpha g (V, Y) + \beta q (V) q(Y)$\\
\vspace{.2cm}
(3.8) \hspace{.5in} $|u(Q)|^2 = - \cfrac{\alpha}{\beta} (1 - \lambda^2)$

\vspace{.2cm}
\parindent=8mm
Also on a cylindrical hypersurface $M$ with $(\phi, g, u, v, \lambda)-$ structure of a Sasakian manifold $\tilde M$, we have\\
\vspace{.2cm}
(3.9) \hspace{1.5in} $v (\phi Q) = 0.$\\
\vspace{.2cm}
(3.10) \hspace{1.45in} $u(Q) = 0 \Leftrightarrow q(U) = 0.$\\\\
\begin{center}
{\bf MAIN RESULT}
\end{center}
 Now we prove\\
{\bf {Theorem~3.1~:}}~ {\it {On the Quasi umbilical hypersurface $M$ with $(\phi, g, u, v, \lambda )-$ structure of a Sasakian manifold $\tilde M$, we have}}\\
\vspace{.2cm}
(a) \hspace{1.5in} $ q(V) = v (Q) = 0.$\\
\vspace{.2cm}
(b) \hspace{1.5in} $w (U) = \cfrac{1 - \lambda^2}{\lambda} - \cfrac{U \lambda}{\lambda}.$\\
\vspace{.2cm}
(c) \hspace{1.5in} $w (V) = \cfrac{\alpha \, (\lambda^2 - 1)}{\lambda} - \cfrac{V \lambda}{\lambda}.$\\
{\it {Proof~}:}~Put $Y = U$ in ((3.1), we have\\
\vspace{.2cm}
\hspace{1.85in}$h (X, U) = \alpha u (X) + \beta q (X) q (U)$\\
\vspace{.2cm}Using (2.18), we get
$$\alpha u (X) + \beta q (X) q (U) = 0.$$
Put $X = U$, then
$$\alpha (1 - \lambda^2) + \beta |q (U) |^2 = 0$$
$$| q (U) |^2 = - \frac{\alpha}{\beta} \, (1 - \lambda^2)$$
If we take $X = V$, then we get
$$\beta q (V) \, q(U) = 0$$
since $\beta \ne 0$, $q(U) \ne 0$ therefore
$$q(V) = 0 = v(Q)$$
which proves (a).\newline
Also from (3.1) and (2.16), we get

\vspace{.2cm}
\noindent
(3.11) \hspace{1in} $\alpha g (Y, V) + \beta q (Y) q (V) = u (Y) - Y \lambda - \lambda w (Y)$.\\
\vspace{.2cm}
Let $Y = U$ and using $g(U, V) = u(V) = 0$, we get
$$w(U) = \frac{1 - \lambda^2}{\lambda} - \frac{U \lambda}{\lambda}$$
which proves (b).\newline
\vspace{.2cm}
Further putting $Y = V$ in (3.11), we get $$w (V) = \cfrac{\alpha \, (\lambda^2 - 1)}{\lambda} - \cfrac{V \lambda}{\lambda}$$
which proves (c).\\
\vspace{.2cm}
If Quasi umbilical hypersurface is cylindrical then from (c), we have $$w = - d (\log \lambda ).$$\\
Therefore, we can state the following theorem as a corollary of Theorem 3.1:\\\\
{\bf{Theorem~3.2~:}}~ {\it {If Quasi umbilical hypersurface $M$ of a Sasakian manifold $\tilde M$ with $(\phi, g, u, v, \lambda)-$structure is cylindrical then $1-form \, w$ satisfies the following relation}}
 $$w = - d (\log \lambda).$$\\

\vspace{.2cm}
\hspace{.5cm}{\it{ Authors' addresses}:}\\
{Sachin Kumar Srivastava \newline
 Department of Mathematics \newline
Central University of Himachal Pradesh \newline
Dharamshala-176215, INDIA. }\\ 
email: sachink.ddumath@gmail.com.\\\\ Alok Kumar Srivastava\\
Department of Mathematics\newline Govt. Degree College \newline Chunar -
Mirzapur,U.P., India. \\
email: aalok\_sri@yahoo.co.in.
\end{document}